\documentstyle[11pt]{article}
\newtheorem{theo}{Theorem}[section]
\newtheorem{definition}{Definition}[section]
\newtheorem{prop}[theo]{Proposition}
\newtheorem{lemma}[theo]{Lemma}
\newtheorem{coro}[theo]{Corollary}

\newcommand{\SS}{{\cal S}}
\newcommand{\SSS}{(S,\overline{S})}

\begin{document}
\date{}

\title{
Testing Equality in Communication Graphs
}

\author{Noga Alon 
\thanks{Sackler School of Mathematics 
and Blavatnik School of
Computer Science, Tel Aviv University, Tel Aviv 69978, Israel.
Email: {\tt nogaa@tau.ac.il}.  Research supported in part by a
USA-Israeli
BSF grant 2012/107, by an ISF grant 620/13 and
by the Israeli I-Core program.
}
\and
Klim Efremenko
\thanks{
Blavatnik School of
Computer Science, Tel Aviv University, Tel Aviv 69978, Israel.
Email: {\tt klimefrem@gmail.com}.  Research supported in part by an
ERC Advanced grant.
}
\and
Benny Sudakov
\thanks{Department of Mathematics, ETH, 8092 Zurich.
Email: {\tt benjamin.sudakov@math.ethz.ch}.
Research supported in part by SNSF grant 200021-149111.
}
}

\maketitle

\begin{abstract}
Let $G=(V,E)$ be a connected undirected graph with $k$ vertices. Suppose
that on each vertex of the graph there is a player having an $n$-bit
string. Each player is allowed to communicate with its neighbors according
to an agreed communication protocol, and the players must decide,
deterministically, if their inputs are all equal. What is the minimum
possible total number of bits transmitted in a protocol solving
this problem ? We determine this minimum up to a lower order
additive term in many cases (but not for all graphs).  
In particular, we show that it is $kn/2+o(n)$ for any
Hamiltonian $k$-vertex graph, and that for any $2$-edge connected 
graph with $m$ edges containing no
two adjacent vertices of degree exceeding $2$ it is $mn/2+o(n)$.
The proofs combine graph theoretic ideas with
tools from additive number theory.
\end{abstract}

\section{The problem}

Let $G=(V,E)$ be a connected undirected graph with $k$ vertices. Suppose
that on each vertex of the graph there is a player having an $n$-bit
string. Each player is allowed to communicate with its neighbors according
to an agreed communication protocol, and the players must decide,
deterministically, whether or not
their inputs are all equal. In a trivial protocol
the players fix a rooted spanning tree of the graph, and each of them,
besides the one at the root, transmits his bits to his parent, and
each one (including the root)
checks 
that his input is equal to those he received from 
each of his children. This
shows that a total communication of roughly $(k-1)n$ bits suffices.
Somewhat surprisingly, it turns out that for complete graphs $G$ with at
least $3$ vertices one can do better. It is shown in \cite{LV} that for
$G=K_k$ at least $kn/2$ bits of communication are needed, and the authors
also obtain a nontrivial upper bound (which is not tight). Brody \cite{Br}
has used the graphs constructed in \cite{AMS} to show that for $G=K_3$,
$3n/2 + o(n)$ bits suffice, showing that the lower bound is tight in this
case up to a low order additive error term. In \cite{AMS} we mentioned
(without giving a detailed proof) that we can use a hypergraph extension
of the construction in \cite{AMS} to show that for $G=K_k$ the minimum
possible number of bits in a communication protocol for the above problem
on $G$ is $(1 + o(1))k n/2.$ Brody and H{\aa}stad have independently found
a similar protocol, using the $k$-cliques of the graphs in \cite{AMS}.

Here we consider the case of general graphs $G$, obtaining upper
and lower bounds which are nearly tight in many (but not all)
cases. Our upper bounds are based on an extension of the graphs of
Rusza and Szemer\'edi \cite{RS}, 
similar to the extension given in \cite{Al}. We also
observe that linear communication protocols cannot improve the
trivial upper bound. Finally, we suggest two
competing conjectures about the possible
answer for every graph.

Let $f(n,G)$ denote the minimum number of bits transmitted in a
communication protocol solving the problem  on $G$. It is clear
that the function $f(n,G)$ is sub-additive, and hence by Fekete's
Lemma the limit of the ratio $f(n,G)/n$ as $n$ tends to infinity
exists. Denote this limit by $f(G)$. The parameter $f(G)$ is the
main object of study in the present short paper.

\section{Results}

Recall that a {\em block} of a graph is a maximal two-connected
subgraph, where every bridge is also a block. 
It is well known that any graph is the edge-disjoint
union of its blocks, and the vertices belonging to more than one
block are the cut vertices of the graph.
Our first observation is the following.
\begin{prop}
\label{p21}
For any connected graph $G$ with blocks $G_1,G_2, \ldots ,G_s$,
$$
f(G)=\sum_{i=1}^s f(G_i).
$$
\end{prop}

For a connected graph $G$ let $c_2(G)$ denote the minimum number of
edges in a $2$-edge connected spanning subgraph of $G$, where the
subgraph is allowed to contain the same edge of $G$ twice. Thus,
for any tree $G$ with $k$ vertices $c_2(G)=2(k-1)$ and for any
graph $G$ with $k$ vertices $c_2(G)=k$ if and only if $G$ is
Hamiltonian. Our main upper bound  for $f(G)$ is the following.
\begin{theo}
\label{t22}
For any connected graph $G$, $f(G) \leq 0.5 c_2(G)$
\end{theo}
\begin{definition}
For a connected graph $G$ let $\SS$ denote the set of all cuts in
$G$. For any edge  $e$ of $G$ let $\SS_e$ denote the set of all
cuts containing $e$.
A fractional packing of cuts in $G$ is a function
$g: \SS \mapsto [0,1]$ so that for every edge $e $ of $G$,
$\sum_{e \in \SS_e}  g(S,\overline{S}) \leq 1$.
Let $fc(G)$ denote the maximum 
possible value of $\sum_{(S,\overline{S}) \in \SS}
g(S,\overline{S})$, where the maximum is taken over all fractional
packings of cuts $g$.
\end{definition}
\begin{theo}
\label{t23}
For any connected graph $G$, $f(G) \geq fc(G)$.
\end{theo}
Note that this implies that $f(G) \geq k/2$ for any $k$-vertex
graph, as the function assigning to all cuts determined by a single
vertex the value $1/2$ is always a fractional packing of cuts.
Note also that clearly if $G'$ is a spanning  subgraph of $G$ then
$f(G') \geq f(G)$ and hence the above $k/2$ lower bound also follows
from the fact that $f(K_k)=k/2.$

By the last theorem $f(G) \geq \alpha(G)$ for every $G$. The 
two theorems above suffice to determine $f(G)$ in many cases.
\begin{coro}
\label{c24} $~~~$
\begin{enumerate}
\item
For any Hamiltonian graph $G$ with  $k$ vertices 
$f(G)=k/2$.
\item
For any complete bipartite graph $G=K_{s,t}$ with $t \geq s \geq
1$, $f(G)=t$.
\item
For any $2$-edge connected graph $G$ in which no two vertices of degree
bigger than $2$ are adjacent, $f(G)$ is exactly half the number of
edges of $G$.
\end{enumerate}
\end{coro} 

A communication protocol is called {\em linear} if any bit it
transmits is a linear combination of the input bits (and the bits
received already). For simplicity we consider only linear
combinations over $Z_2$, but the (simple) result that follows 
can be easily extended to all finite fields.
\begin{prop}
\label{p25}
For any connected graph $G$ on $k$ vertices, any linear protocol
for solving the equality problem requires communication of at least
$(k-1)n$ bits.
\end{prop}

\section{Proofs}

\subsection{Preliminaries}

We start with the simple proofs of Propositions \ref{p21}
and \ref{p25}

\vspace{0.2cm}

\noindent
{\bf Proof of Proposition \ref{p21}:}\,
We apply induction on the number of blocks $s$. For $s=1$ there is
nothing to prove. Assuming the result holds for $s-1$ we prove it
for $s$, $s \geq 2$. Let $G, G_1, \ldots ,G_s$ be as in the
proposition, and assume, without loss of generality, that $G_s$ is
an end-block.  Let $v$ be the unique cut-vertex in $G_s$ and let
$G'$ be the graph obtained  from $G$ by removing all vertices of
$G_s$ besides $v$. Thus $G'$ has $s-1$ blocks $G_1,G_2, \ldots
,G_{s-1}$.

To show that $f(G) \leq \sum_{i=1}^s f(G_i)$ observe that one can
first apply the best protocol for solving the problem in $G_s$. If 
all vertices of $G_s$ have the same bit string as $v$, we can now
apply the best protocol for $G'$ to complete the required task,
thus establishing the upper bound.

To prove the lower bound consider the best protocol for solving the
problem for $G$. By considering its behavior only on inputs of
length $n$ in
which all vertices of $G_s$ have equal inputs we conclude that the
number of bits transmitted by this protocol along edges of $G'$ is
at least $f(n,G')$. Similarly, by considering the scenarios in
which all vertices of $G'$ have the same strings we conclude that
the number of bits transmitted along edges of $G_s$ is at least
$f(n,G_s)$. This establishes the lower bound, 
completing the proof. \hfill $\Box$
\vspace{0.2cm}

\noindent
{\bf Proof of Proposition \ref{p25}:}\,
Consider a linear protocol for the problem, and suppose it
transmits $m$ bits. Each bit is a linear combination of the $nk$
bits representing the inputs of the $k$ vertices. For each such
combination, define a linear equation equating it to zero. The set
of all these $m$ equations is a homogeneous system of $m$ linear
equations in $kn$ variables. If $m<(k-1)n$ then the dimension of
the solution space is bigger than $n$. However, the dimension of
the space of all inputs in which  all strings are equal is 
$n$, hence there is a solution, call it $s$,
 in which not all input strings are
equal. Note that if each input  string is the $0$ vector, then all
bits transmitted are $0$, and the protocol must accept. Therefore,
it must also accept the input $s$, as with this input all bits transmitted
are also zero. But this means that the protocol errs on the input
$s$, showing that a total communication of less than
$(k-1)n$ is impossible in the linear case, as needed.  \hfill
$\Box$

\subsection{The upper bound}

In this section we prove Theorem \ref{t22}. We need several lemmas,
the first one is a known extension of the construction of Behrend
in \cite{Be} of dense sets of integers with no $3$-term arithmetic
progressions. 

A linear equation with integer coefficients

\begin{equation}
\label{e30}
\sum a_{i} x_{i} = 0
\end{equation}

\noindent
in the unknowns $x_{i}$ is {\em homogeneous} if $\sum a_{i} =
0$.
If $X\subseteq M= \{1,2,\ldots,m\}$, we say that $X$
{\em has no non-trivial
solution to} (\ref{e30}), if whenever $x_{i}\in X$
and $\sum a_{i}x_{i} = 0$, it follows that all $x_{i}$ are equal.
Thus, for example, $X$ has no nontrivial
solution to the equation $x_1-2x_2+x_3=0$ iff it contains no
three-term arithmetic progression.
\begin{lemma}[see, e.g., \cite{Al}, Lemma 3.1]
\label{l31}
For every fixed integer $k \geq 2$   and every positive integer
$m$,
there exists
a subset $X \subset M=\{1,2, \ldots ,m\}$ of size at least
$$
|X| \geq \frac{m}{e^{10 \sqrt {\log m \log k}}}
$$
with no non-trivial solution to the equation
\begin{equation}
\label{e31}
x_1+x_2 + \ldots +x_k=k x_{k+1}.
\end{equation}
\end{lemma}
Note that if there is no nontrivial solution for the above equation
there is also no non-trivial solution for each of the equations
$x_1+x_2 + \ldots +x_r=r x_{r+1} $ for $r \leq k$, since a
non-trivial solution of that together with
$x_{r+1}=x_{r+2}  = \ldots = x_k=x_{k+1}$ yields 
a non-trivial solution of
(\ref{e31}).

We also need a basic result on $2$ connected graphs, first proved
by Whitney \cite{Wh}. 
An {\em ear} of an 
undirected graph $G$ is a path $P$ where the two endpoints
of the path may coincide, but where otherwise no repetition of
edges or vertices is allowed.
A {\em proper ear decomposition} of $G$ is a
partition of its set of edges into a sequence of ears, such that
the first ear is a cycle, the two endpoints of any other ear are
distinct
and belong to earlier ears in the
sequence and the internal vertices of each ear (if any) do not
belong to any earlier ear. The following result was first proved by
Whitney (it is also an easy consequence of Menger's Theorem.)
\begin{lemma}[Whitney \cite{Wh}]
\label{l32}
A graph $G$ 
is $2$ connected if and only if it has a proper ear
decomposition.
\end{lemma}

Let $H$ be a graph with $k$ vertices $\{v_1,v_2, \ldots ,v_k\}$.
Let $F$ be  a $k$-partite graph with classes of vertices $V_1,V_2, \ldots
,V_k$. A copy of $H$ in $F$ is called a {\em special copy} if for
each $1 \leq i \leq k$ the vertex playing the role of $v_i$ belongs
to $V_i$. Call $F$ a 
{\em faithful host} for $H$ if the set of its edges is the edge-disjoint
union of special copies of $H$, and 
$F$ contains no other special copy of $H$ besides the
$|E(F)|/|E(H)|$ copies defining its set of edges.
The following lemma is a crucial ingredient in the proof of 
Theorem  \ref{t22}. The special case when $H$ is a cycle is proved
in \cite{Al}.
\begin{lemma}
\label{l33}
Let $H$ be a $2$-connected graph with $k$ vertices, and let $m$ be
a positive integer. Then there is 
a faithful host $F$ for $H$ with classes of vertices $V_1, \ldots
,V_k$, each of size $km$, containing at least
$$
\frac{m^2}{e^{10 \sqrt {\log m \log k}}}
$$
special copies of $H$.
\end{lemma}
{\bf Proof:}\, 
By Lemma \ref{l32} there is a proper ear decomposition of $H$. Fix
such a decomposition, and denote the ears in it by $P_1,P_2,
\ldots, P_s$, in order, where $P_1$ is a cycle and each $P_j$ for
$j>1$ is a path whose endpoints lie on vertices of earlier ears.
Define a numbering of the vertices of $H$ as follows. The vertices
of the cycle $P_1$ are numbered $v_1,v_2, \ldots ,v_t$, according
to their order on the cycle. Assuming we have already numbered all
vertices in the first $p$ ears by $v_1, v_2, v_3 \ldots ,v_{\ell}$,
consider the next ear $P_{p+1}$. If it contains no internal
vertices there is no new vertex in it that should be numbered.
Otherwise, suppose the endpoints of this ear are $v_i$ and $v_j$,
where $i<j$, and suppose it has $q$ internal vertices. Then this
ear is a path of length $q+1$ from $v_i$ to $v_j$ and its vertices
are numbered so that the vertices of the path are 
$v_j,v_{\ell+1},v_{\ell+2}, \ldots ,v_{\ell+q}, v_i$ in this order.

Let $X \subset \{1,2, \ldots ,m\}$ be as in Lemma \ref{l31}. The
host graph $F$ is defined as follows. Its vertex classes are the
classes $V_1,V_2, \ldots ,V_k$, where each $V_i$ is of size $km$
(the first classes can be smaller, but this is not essential for our
purpose here, hence we prefer the more symmetric description as
above). With slight abuse of notation  denote the vertices of each
set $V_i$ by $\{1,2, \ldots ,km\}$ but recall that these sets are
pairwise disjoint. The graph $F$ contains $m|X|$ special copies of
$H$ defined as follows. For each integer $y$, 
$1 \leq y \leq m$ and each $x \in X$, there is a special 
copy of $H$ in $F$, which we denote by $H_{x,y}$,
in which $y+(i-1)x \in V_i$ is 
the vertex playing the role of $v_i$ (for all $1 \leq i \leq k$).
It is easy to see that all these special copies are pairwise edge  
disjoint. In fact, these copies satisfy a stronger property: no
two of them share two vertices, since the values of $y+(i-1)x$ for
two distinct indices $i$ determine uniquely $x$ and $y$. It remains
to prove that the only special copies of $H$ in $F$ are the copies
$H_{x,y}$ used in its definition. Let $H'$ be such a special copy.
Then it contains an edge between $V_1$ and $V_2$ which connects
$y \in V_1$ to $y+x \in V_2$, where $1 \leq y \leq m$ and $x \in
X$. Let $u_1, u_2, \ldots ,u_k$ be the vertices of $H'$, where
$u_i \in V_i$ for all $i$. Our objective is to prove that
$u_i =y+(i-1)x$ for all $i$. To do so we show, by induction on $p$,
that this holds for each of the vertices $u_i$ where $v_i$ 
belongs to the
union of the vertices in the first $p$ ears in the ear
decomposition of $H$. The first ear, $P_1$, is a cycle
on the vertices $v_1,v_2,\ldots ,v_t$. By the definition of $F$
there are $x_1=x,x_2,\ldots ,x_t \in X$ so that 
$u_{i+1}-u_i=x_i$ for all $1 \leq i \leq t-1$ and $u_t-u_1=
(t-1)x_t$. Therefore
$x_1+x_2+ \ldots +x_{t-1}=(t-1)x_t$. Since $t \leq k$, the property
of the set $X$ implies that $x_i=x_1=x$ for all $1 \leq i \leq t$,
establishing the required beginning of the induction. Assuming the 
induction claim holds for the vertices in the first $p$ ears,
consider the next ear $P_{p+1}$. If it contains no internal
vertices there is nothing to prove, hence assume it contains
$q$ internal vertices. Let the ear $P_{p+1}$ be
$v_j,v_{\ell+1},v_{\ell+2}, \ldots ,v_{\ell+q},v_i$, where
$i <j$. By the induction hypothesis 
$u_i=y+(i-1)x$ and $u_j=y+(j-1)x$. By the construction of $F$
there are $x_1,x_2, \ldots ,x_{q+1} \in X$ so that
$u_{\ell+1}-u_j=(\ell+1-j)x_1$, $u_{\ell+i+1}-u_{\ell+i}=x_{i+1}$
for $1 \leq i \leq q-1$, and $u_{\ell+q}-u_i=(\ell+q-i)x_{q+1}$.
Since
$$
(u_j-u_i)+(u_{\ell+1}-u_j)+(u_{\ell+2}-u_{\ell+1}) +
\ldots +(u_{\ell+q}-u_{\ell+q-1})=u_{\ell+q}-u_i
$$
we conclude that
$$
(j-i)x+(\ell+1-j)x_1+x_2+ \ldots +x_q=(\ell+q-i)x_{q+1}.
$$
As $\ell+q-i \leq k$ the property of $X$ implies 
that $x=x_1=x_2=\ldots =x_{q+1}$ completing the proof of the
induction and implying the assertion of the lemma. \hfill 
$\Box$
\vspace{0.2cm}

\noindent
{\bf Proof of Theorem \ref{t22}:}\, 
Let $G'$ be a two edge-connected spanning subgraph of  $G$
with $c_2(G)$ edges. It may 
contain two copies of some edges, but by the minimality in the
definition of $G'$ this 
is the case only for bridges of (the underlying subgraph of) 
$G'$. We have to show that
$f(G')$ is at most half the number of its edges. By Proposition
$\ref{p21}$ it suffices to prove it for all blocks of 
$G'$, where for blocks consisting of a single edge (taken twice)
this is trivial, as obviously $f(K_2)=1$. 
Every nontrivial block of $G'$ is $2$ connected, and it thus
suffices to show that for any $2$-connected graph $H=(V,E)$,
$f(H) \leq 0.5 |E|$.  Let $k$ denote the number of vertices of $H$.
For a given (large) integer $n$, let  $m$ be the smallest integer 
so that 
$$
\frac{m^2}{e^{10 \sqrt {\log m \log k}}} \geq 2^n.
$$
Thus 
$$ 
\log_2 m=0.5 n +O(\sqrt {n  \log k})
$$
and 
$$
\lceil \log_2 (k m)  \rceil =0.5 n+O(\sqrt {n  \log k})+O(\log k)
=(0.5+o(1)) n.
$$
Fix a numbering $v_1,v_2, \ldots ,v_k$ of the 
vertices of $H$ according to the proof of Lemma \ref{l33}, and let
$F$ be a faithful host for $H$, with classes of vertices 
$V_1,V_2, \ldots ,V_k$, containing at least $2^n$
special copies of $H$. 
Fix $2^n$ special copies. The input strings are now represented by
special copies of $H$ in $F$. Orient the edges of $H$ so that the
indegree of every vertex is positive. This is possible, since 
$H$ is $2$-connected. Indeed, using an ear decomposition of $H$
we can orient the initial cycle cyclically and then orient each ear
as a directed path.
The player $P_i$ residing at the
vertex $v_i$ of $H$ transmits the identity  of the vertex $u_i$ in
the special copy of $H$ representing his input to all players
$P_j$ so that there is an edge of $H$ oriented from $v_i$ to $v_j$.
Note that this amounts  to a total transmission of 
$$
\lceil \log_2 (k m)  \rceil |E(H)|=(0.5+o(1))|E|n
$$
bits.
In addition, 
each player observes if the identities of the vertices he 
received from his inneighbors are indeed consistent with the ones
in 
his copy, and reports about this to his outneighbors (this amounts
to another single bit per edge). If there is some inconsistency,
this information reaches some player who reports that the
inputs are not all equal. If everything is consistent, the players
report that all inputs are equal.

It is clear that if all inputs are equal then the players report
so. To complete the proof we show that if they report that the
inputs are all equal, this is indeed the case. For every
$i$ let 
$u_i$ be the identity of the vertex in $V_i$ reported 
by $i$ to his outneighbors. Let the
special copies of the players be $H_1,H_2, \ldots ,H_k$, where
$H_i$ is the copy of the player  $P_i$. If $(v_j,v_i)$ is an edge
of $H$ oriented from $v_j$ to $v_i$, and $v_i$ who gets the
identity of the vertex $u_j \in V_j$ from the player $P_j$, finds it
consistent with his copy, then the edge $u_ju_i$ belongs to the
special copy $H_i$ of $P_i$. Therefore, if no player reports an
inconsistency, then the subgraph of $F$ on the vertices $u_1,u_2,
\ldots ,u_k$ is a special  copy of $H$ in $F$. However, since
$F$ is a faithful host for $H$ this copy must be one of the
original special copies of $H$ in $F$, and as it contains an edge
of each $H_i$ (as the indegree of each vertex is positive) 
this special copy must be equal to $H_i$ for all $i$, showing that
indeed all these copies are equal. This completes the proof. \hfill
$\Box$

\subsection{The lower bound}
\vspace{0.2cm}

\noindent
{\bf Proof of Theorem \ref{t23}:}\, 
Consider a deterministic communication protocol that solves the
equality  problem for inputs with $n$ bits 
on $G=(V,E)$. For each edge $e \in E$, let
$b(e)$  denote the number of bits transmitted during the protocol
along $e$. We claim that for every cut  $(S,\overline{S})$ in $G$
$\sum_{e \in \SSS} b(e) \geq n$. Indeed, otherwise there are two
distinct strings of length $n$, $x$ and $y$, so that the
communication along the edges of the cut is identical when all inputs
are $x$ and when all  inputs are $y$. But in that case it is easy
to see that the protocol behaves identically when
all inputs are $x$, when all inputs are $y$, and also when all
vertices of $S$ have input $x$ and all those in $y$ have input $y$
(and vice versa). Thus the protocol cannot behave correctly,
proving the claim.

By the claim it follows that  a lower bound for $f(n,G)$ is 
the solution of the following linear program:
\begin{equation}
\label{e32}
\mbox{Minimize} ~~\sum_e b(e)~~ \mbox{subject to the constraints}
\end{equation}
$$
b(e) \geq 0~~ \mbox{for all}~~ e \in E~~ \mbox{and}
$$
$$
\sum_{e \in \SSS} b(e) \geq n~~ \mbox{for every cut}~~ \SSS \in \SS,
$$
where $\SS$ is the set of all cuts of $G$.

The dual of this program is:
\vspace{0.2cm}

\noindent
Maximize $n \cdot \sum_{\SSS \in \SS} g\SSS$ subject to the constraints
\vspace{0.1cm}

\noindent
$g\SSS \geq 0$ for all $\SSS \in \SS$ and 
$\sum_{\SSS, e \in \SSS} g\SSS \leq 1$ for every edge $e \in E$.
\vspace{0.1cm}

\noindent
This last maximum is exactly $n \cdot fc(G)$, 
completing the proof. \hfill $\Box$
\vspace{0.2cm}

\noindent
{\bf Proof of Corollary \ref{c24}:}\, 
\begin{enumerate}
\item
By Theorem \ref{t23} and the paragraph following its statement
$f(G) \geq k/2$ for any $k$-vertex graph $G$.
By Theorem \ref{t22},
for the cycle $C_k$ on $k$ vertices $f(C_k) \leq k/2$. The desired
result follows since if $G'$ is a spanning subgraph of $G$ then
clearly $f(G) \leq f(G')$.
\item
For any tree $T$ on $k$  vertices $f(T)=k-1$ (for example, by Proposition 
\ref{p21}). This implies the result for $s=1$. For larger $s$ the 
lower bound follows from Theorem \ref{t23}
by the fact that for any graph $G$ with
independence number $\alpha=\alpha(G)$, $fc(G) \geq \alpha$ as
the $\alpha$ cuts $(v, V(G)-\{v\})$ for $v$ in a maximum
independent set are pairwise edge disjoint. The upper bound follows
from Theorem \ref{t22} by considering a spanning subgraph of
$K_{s,t}$ consisting of a cycle of length $2s$  together with two
of the edges incident with any vertex of $K_{s,t}$  uncovered by
the cycle.
\item
The upper bound follows from Theorem \ref{t22}. To prove the lower
bound note that $G$ is the edge disjoint union of induced 
paths, each of length at least $2$. For each such path 
$v_1,v_2,  \ldots ,v_s$ in which all internal vertices are of
degree $2$ in $G$, consider the cuts $(v_i,V-\{v_i\})$ for all
$1 < i <s$, and the cut 
$$
(\{v_2,v_3 \ldots ,v_{s-1}\},V-\{v_2,v_3 \ldots ,v_{s-1}\})
$$
(if $s=3$ we take the same cut twice). This is a collection of
$|E(G)|$ cuts covering each edge  exactly twice, hence $fc(G)
\geq |E(G)|/2$, as shown by giving each of these cuts weight $1/2$.
This completes the proof. \hfill $\Box$
\end{enumerate}

\section{Open problems}

\begin{itemize}
\item
Is $f(G)=0.5 c_2(G)$ for any connected graph $G$ ? 
\item
If not, is $f(G)=fc(G)$ for any connected graph $G$ ? 
\item
Is it true that
for a graph $G$ on $k$ vertices $f(G)=k/2$ if and only if $G$ is
Hamiltonian? (Note that if this is the case, then the computational
problem of computing $f(G)$ for a given input  graph $G$ is
NP-hard.)
\item
It is not difficult to show that for any $d$-regular graph $G$
on $k$ vertices
which is also $d$-edge connected, $fc(G)=k/2$. Indeed, as mentioned 
in the paragraph  following the statement of Theorem \ref{t23},
$fc(G) \geq k/2$ for any $k$ vertex graph. To prove the upper bound
note that for any $d$-regular $d$ edge-connected graph $G=(V,E)$, the
function $b(e)=n/d$ for every edge $e \in E$ is a solution of the
linear program (\ref{e32}). 

Thus, for any such $G$ the lower bound
for $f(G)$ provided by Theorem \ref{t23} is $k/2$ whereas if it is
not Hamiltonian the upper bound provided by Theorem \ref{t22} is
strictly larger.

A specific interesting example is
the Petersen graph $P$ which is $3$-regular, $3$-connected and
non-Hamiltonian. Indeed 
$c_2(P)=11$ and $fc(G)=5$, implying that
$$
5 \leq f(P) \leq 5.5
$$ 
What is $f(P)$ ?
\end{itemize}

\end{document}